\documentclass[11pt]{article}
\usepackage{amscd,amsmath,amsthm,amssymb,xypic}

\theoremstyle{plain}  

\newtheorem{thm}{Theorem}[section]
\newtheorem{prop}[thm]{Proposition}
\newtheorem{lem}[thm]{Lemma}   
\newtheorem{cor}[thm]{Corollary}

\theoremstyle{definition}

\newtheorem{prop-defn}[thm]{Proposition--Definition}

\theoremstyle{remark}

\newtheorem{rem}[thm]{Remark}

\DeclareMathOperator{\cHom}{\mathcal{H}\mathnormal{om}}
\DeclareMathOperator{\Ext}{Ext}
\DeclareMathOperator{\cExt}{\mathcal{E}\mathnormal{xt}}

\DeclareMathOperator{\Pic}{Pic}

\DeclareMathOperator{\Sing}{Sing}

\DeclareMathOperator{\pr}{pr}

\newcommand{\QED}{\ifhmode\unskip\nobreak\fi\quad {\rm Q.E.D.}} 

\newcommand{\bA}{\mathbb A}
\newcommand{\bC}{\mathbb C}

\newcommand{\bN}{\mathbb N}
\newcommand{\bP}{\mathbb P}
\newcommand{\bQ}{\mathbb Q}

\newcommand{\cB}{\mathcal B}

\newcommand{\cX}{\mathcal X}

\newcommand{\sD}{\mathcal D}

\newcommand{\cF}{\mathcal F}

\newcommand{\cO}{\mathcal O}
\newcommand{\sL}{\mathcal L}
\newcommand{\cM}{\mathcal M}
\newcommand{\cN}{\mathcal N}

\newcommand{\cU}{\mathcal U}

\newcommand{\omgn}{\overline{\mathcal M}_{g,n}}
\newcommand{\oMgn}{\overline{M}_{g,n}}
\newcommand{\om}{\overline{\mathcal M}}
\newcommand{\oM}{\overline{M}}
\newcommand{\ocM}{\overline{\mathcal M}}

\title{The moduli space of curves is rigid}
\author{Paul Hacking}
\begin{document}
\maketitle

\begin{abstract}
We prove that the moduli stack $\omgn$ of stable curves of genus $g$ with $n$ marked points is rigid, that is, has no infinitesimal deformations. This confirms the first case of a principle proposed by Kapranov. It can also be viewed as a version of Mostow rigidity for the mapping class group.
\end{abstract}

\section{Introduction}

Kapranov has proposed the following informal statement \cite{Kapranov97}.
Given a smooth variety $X=X(0)$, 
consider the moduli space $X(1)$ of varieties obtained as deformations of 
$X(0)$, the moduli space $X(2)$ of deformations of $X(1)$, and so on. 
Then this 
process should stop after $n = \dim X$ steps, i.e., $X(n)$ should be rigid (no 
infinitesimal deformations).
Roughly speaking, one thinks of $X(1)$ as $H^1$ of a sheaf of non-abelian groups on $X(0)$.
Indeed, at least the tangent space to $X(1)$ at $[X]$ is identified with $H^1(T_X)$, where $T_X$ 
is the tangent sheaf, the sheaf of first order infinitesimal automorphisms of $X$.
Then one regards $X(m)$ as a kind of non-Abelian $H^m$, and the analogy with the usual
definition of Abelian $H^m$ suggests the statement above.

In particular, the moduli space of curves should be rigid. In this paper, we 
verify this in the following precise form: the moduli stack of stable 
curves of genus $g$ with $n$ marked points is rigid for each $g$ and $n$.

On the other hand, moduli spaces of surfaces should have non-trivial 
deformations in general. A simple example (for surfaces with boundary) is given in 
Sec.~\ref{example}. It seems plausible that there should be a non-trivial deformation of
a moduli space of surfaces whose fibres parametrise ``generalised surfaces'' in some sense, 
for example non-commutative surfaces. From this point of view the result of this paper 
says that the concept of a curve cannot be deformed.

Let us also note that our result can be thought of as a version of Mostow rigidity for the mapping class group.
Recall that the moduli space $M_g$ of smooth complex curves of genus $g$ is the 
quotient of the Teichm\"{u}ller space $T_g$ by the mapping class group $\Gamma_g$.
The space $T_g$ is a bounded domain in $\bC^{3g-3}$, which is homeomorphic to a ball,
and $\Gamma_g$ acts discontinuously on $T_g$ with finite stabilisers. 
We thus obtain $M_g$ as a complex orbifold with orbifold fundamental group $\Gamma_g$.
The space $T_g$ admits a natural metric, the Weil--Petersson metric, which has negative 
holomorphic sectional curvatures. So, roughly speaking, $M_g$ looks like a quotient of 
a complex ball by a discrete group $\Gamma$ of isometries, with finite volume. Mostow rigidity 
predicts that such a quotient is uniquely determined by the group $\Gamma$ up to complex conjugation.
(This is certainly true if $\Gamma$ acts freely with compact quotient, see \cite{Siu80}.)
In particular, it should have no infinitesimal deformations. Unfortunately I do not know a proof along these lines.

\medskip
\noindent
\emph{Acknowledgements}: We thank M.~Kapranov for suggesting the problem and subsequent encouragement.
We also thank I.~Dolgachev, G.~Farkas, S.~Grushevsky, S.~Keel, J.~Koll\'{a}r, S.~Kov\'{a}cs, R.~Lazarsfeld, A.~Marian, M.~Olsson, 
S.P.~Smith and T.~Stafford for helpful discussions. The author was partially supported by NSF grant DMS-0650052.

\section{Statements}

We work over an algebraically closed field $k$ of characteristic zero.
Let $g$  and $n$ be non-negative integers such that $2g-2+n > 0$.
Let $\omgn$ denote the moduli stack of stable curves of genus $g$ with $n$ marked points.
The stack $\omgn$ is a smooth proper Deligne--Mumford stack of dimension $3g-3+n$.

\begin{thm}\label{mainthm}
The stack $\omgn$ is rigid, that is, has no infinitesimal deformations.
\end{thm}

Let $\partial\omgn \subset \omgn$ denote the boundary of the moduli stack, 
that is, the complement of the locus of smooth curves (with its reduced structure).
The locus $\partial\omgn$ is a normal crossing divisor in $\omgn$.
 
\begin{thm}\label{logdefns}
The pair $(\omgn,\partial\omgn)$ has no locally trivial deformations.
\end{thm}

Let $\oMgn$ denote the coarse moduli space of the stack $\omgn$.
The space $\oMgn$ is a projective variety with quotient singularities.

\begin{thm} \label{coarsedefns}
The variety $\oMgn$ has no locally trivial deformations if $(g,n) \neq (1,2),(2,0),(2,1),(3,0)$.
\end{thm}

\begin{rem}
In the exceptional cases, the projection $\omgn \rightarrow \oMgn$ is
ramified in codimension one over the interior of $\oMgn$, 
and an additional calculation is needed to relate the deformations of 
the stack and the deformations of the coarse moduli space (cf. Prop.~\ref{stack-cmod}).
Presumably the result still holds.
\end{rem}

\section{Proof of Theorem~\ref{logdefns}}
Write $\cB$ for the boundary of $\omgn$.
Let $\Omega_{\omgn}(\log \cB)$ denote the sheaf of $1$-forms on $\omgn$ with logarithmic poles along the boundary,
and $T_{\omgn}(-\log \cB)$ the dual of $\Omega_{\omgn}(\log \cB)$.
The sheaf $T_{\omgn}(-\log \cB)$ is the subsheaf of the tangent sheaf $T_{\omgn}$ consisting of vector fields on $\omgn$
which are tangent to the boundary.
In other words, it is the sheaf of first order infinitesimal automorphisms of the pair $(\omgn,\cB)$.
Hence the first order locally trivial deformations of the pair $(\omgn,\cB)$ are identified with the space
$H^1(T_{\omgn}(-\log \cB))$. To prove Thm.~\ref{logdefns}, we show $H^1(T_{\omgn}(-\log \cB))=0$.

Let $\pi : \cU_{g,n} \rightarrow \omgn$ denote the universal family over $\omgn$.
That is, $\cU_{g,n}$ is the stack of $n$-pointed stable curves of genus $g$ together with an extra section (with no smoothness condition).
Let $\Sigma$ denote the union of the $n$ tautological sections of $\pi$.
We define the boundary $\cB_U$ of $\cU_{g,n}$ as the union of $\pi^*\cB$ and $\Sigma$.

Let $\nu \colon \cB^{\nu} \rightarrow \cB$ be the normalisation of the boundary $\cB$ of $\omgn$,
and $\cN$ the normal bundle of the map $\cB^{\nu} \rightarrow \omgn$.
Then we have an exact sequence
$$0 \rightarrow T_{\omgn}(-\log \cB) \rightarrow T_{\omgn} \rightarrow \nu_*\cN \rightarrow 0.$$
Let $\omega_{\pi}$ denote the relative dualising sheaf of the morphism $\pi$. 

\begin{lem} There is a natural isomorphism 
$$\delta : T_{\omgn}(-\log \cB) \stackrel{\sim}{\longrightarrow} R^1\pi_*(\omega_{\pi}(\Sigma)^{\vee}).$$
\end{lem}
\begin{proof}
For a pointed stable curve $(C,\Sigma_C=x_1+\cdots+x_n)$, the space of first order deformations is equal to $\Ext^1(\Omega_C(\Sigma_C),\cO_C)$.
See \cite[p.79--82]{DM69}. The surjection
$$\Ext^1(\Omega_C(\Sigma_C),\cO_C) \rightarrow H^0(\cExt^1(\Omega_C(\Sigma_C),\cO_C)) = \bigoplus_{q \in \Sing C} \cExt^1(\Omega_C(\Sigma_C),\cO_C)_q$$
sends a global deformation of $(C,\Sigma_C)$ to the induced deformations of the nodes.
\'Etale locally at the point $[(C,\Sigma_C)] \in \omgn$, the boundary $\cB$ is a normal crossing divisor with components $B_q$ indexed by the nodes $q$ of $C$ (the divisor $B_q$ is the locus where the node $q$ is \emph{not} smoothed).
The Kodaira--Spencer map identifies the fibre of the normal bundle of $B_q$ at $[(C,\Sigma_C)]$ 
with the stalk of $\cExt^1(\Omega_C(\Sigma_C),\cO_C)$ at $q$. 

We now work globally over $\omgn$. We omit the subscripts $g,n$ for clarity. Consider the exact sequence
\begin{equation}\label{Omega}
0 \rightarrow \pi^*\Omega_{\om} \rightarrow \Omega_{\cU}(\log \Sigma) \rightarrow \Omega_{\cU/\om}(\Sigma) \rightarrow 0.
\end{equation}
For a sheaf $\cF$ on $\cU$, let $\cExt^i_{\pi}(\cF,\cdot)$ denote the $i$th right derived functor of $\pi_* \circ \cHom(\cF,\cdot)$.
Applying $\pi_* \circ \cHom(\cdot,\cO_{\cU})$ to the exact sequence $(\ref{Omega})$, we obtain a long exact sequence with connecting homomorphism
$$\rho \colon T_{\om} \rightarrow \cExt^1_{\pi}(\Omega_{\cU/\om}(\Sigma),\cO_{\cU}).$$
The map $\rho$ is the Kodaira--Spencer map for the universal family over $\om$ and thus is an isomorphism. 
(Note that, for a point $p=[(C,\Sigma_C)] \in \om$, the base change map 
$$\cExt^1_{\pi}(\Omega_{\cU/\om}(\Sigma),\cO_{\cU}) \otimes k(p) \rightarrow \cExt^1(\Omega_C(\Sigma_C),\cO_C)$$
is an isomorphism. Indeed, by relative duality \cite[Thm.~21]{Kleiman80}, it suffices to show that $\pi_*(\Omega_{\cU/\om}(\Sigma) \otimes \omega_{\pi})$
commutes with base change. This follows from cohomology and base change.)

Consider the two exact sequences
$$0  \rightarrow  T_{\om}(-\log \cB)  \rightarrow  T_{\om}  \rightarrow  \nu_*\cN  \rightarrow  0$$
and
$$0 \rightarrow R^1\pi_*(\Omega_{\cU/\om}(\Sigma)^{\vee}) \rightarrow \cExt^1_{\pi}(\Omega_{\cU/\om}(\Sigma),\cO_{\cU}) \rightarrow \pi_*\cExt^1(\Omega_{\cU/\om}(\Sigma),\cO_{\cU}) \rightarrow 0$$
The Kodaira--Spencer map $\rho$ identifies the middle terms, and induces an identification of the right end terms determined by the deformations of the singularities of the fibres of $\pi$.
We thus obtain a natural isomorphism $\delta$ of the left end terms.
Finally, note that $\Omega_{\cU/\om}(\Sigma)^{\vee}=\omega_{\pi}(\Sigma)^{\vee}$ because $\omega_{\pi}(\Sigma)$ is invertible and agrees with 
$\Omega_{\cU/\om}(\Sigma)$ in codimension $1$. This completes the proof.
\end{proof}

The line bundle $\omega_{\pi}(\Sigma)$ is ample on fibres of $\pi$. Hence $\pi_*(\omega_{\pi}(\Sigma)^{\vee})=0$.
Also $R^i\pi_*(\omega_{\pi}(\Sigma)^{\vee})=0$ for $i>1$ by dimensions.
So $H^{i+1}(\omega_{\pi}(\Sigma)^{\vee}) = H^i(R^1\pi_*(\omega_{\pi}(\Sigma)^{\vee}))$ for all $i$ by the Leray spectral sequence.
Hence the isomorphism $\delta$ induces an isomorphism 
\begin{equation} \label{red2U}
H^i(T_{\omgn}(-\log \cB)) \stackrel{\sim}{\longrightarrow} H^{i+1}(\omega_{\pi}(\Sigma)^{\vee})
\end{equation}
for each $i$.

Let $U_{g,n}$ denote the coarse moduli space of the stack $\cU_{g,n}$ and $p : \cU_{g,n} \rightarrow U_{g,n}$ the projection.
The line bundle $\omega_{\pi}(\Sigma)$ on the stack $\cU_{g,n}$ defines a $\bQ$-line bundle $p_*^{\bQ}\omega_{\pi}(\Sigma)$
on the coarse moduli space $U_{g,n}$ (see Sec.~\ref{app}).
We use the following important result, which is essentially due to Arakelov \cite[Prop.~3.2, p.~1297]{Arakelov71}.
We refer to \cite[Sec.~4]{Keel99} for the proof.
\begin{thm} \label{bignef}
The $\bQ$-line bundle $p_*^{\bQ}\omega_{\pi}(\Sigma)$ is big and nef on $U_{g,n}$.
\end{thm} 

It follows by Kodaira vanishing (see Thm~\ref{KVstack}) that 
$H^i(\omega_{\pi}(\Sigma)^{\vee})=0$
for $i < \dim \cU_{g,n}$.
Combining with (\ref{red2U}), we deduce
\begin{prop} \label{logvan}
$H^i(T_{\omgn}(- \log \cB))=0$ for $i < \dim \omgn$.
\end{prop}
In particular, $H^1(T_{\omgn}(- \log \cB))=0$ if $\dim \omgn > 1$. The remaining cases are easy to check.
This completes the proof of Theorem~\ref{logdefns}.

\section{Proof of Theorem~\ref{mainthm}} \label{defstack}

We now prove that $\omgn$ is rigid. 
Since $\omgn$ is a smooth Deligne--Mumford stack, its first order infinitesimal deformations are identified with the space $H^1(T_{\omgn})$, 
and we must show that $H^1(T_{\omgn})=0$.
Consider the  exact sequence 
$$0 \rightarrow T_{\omgn}(-\log \cB) \rightarrow T_{\omgn} \rightarrow \nu_*\cN \rightarrow 0$$
and the associated long exact sequence of cohomology
$$ \cdots \rightarrow H^i(T_{\omgn}(-\log \cB)) \rightarrow H^i(T_{\omgn}) \rightarrow H^i(\cN) \rightarrow \cdots $$
We prove below that $H^i(\cN)=0$ for $i < \dim \cB$. Now $H^i(T_{\omgn}(-\log \cB))=0$ for $i< \dim \omgn$ by Prop.~\ref{logvan},
so we deduce 
\begin{prop} \label{van}
$H^i(T_{\omgn})=0$ for $i< \dim \omgn-1$. 
\end{prop}
In particular, $H^1(T_{\omgn})=0$ if $\dim \omgn > 2$.
In the remaining cases it is easy to check that $H^1(\cN)=0$, so again $H^1(T_{\omgn})=0$.

The irreducible components of the normalisation $\cB^{\nu}$ of the boundary $\cB$ of $\omgn$ are finite images of the following stacks
\cite[Def.~3.8, Cor.~3.9]{Knudsen83a}:
\begin{enumerate}
\item $\om_{g_1,S_1 \cup \{n+1\}} \times \om_{g_2,S_2 \cup \{n+2\}}$ where $g_1+g_2=g$ and $S_1,S_2$ is a partition of $\{1,\ldots, n\}$.
\item $\om_{g-1,n+2}$
\end{enumerate}
Here $\om_{h,S}$ denotes the moduli stack of stable curves of genus $h$ with marked points labelled by a finite set $S$.
In each case the map to $\cB^{\nu}$ is given by identifying the points labelled by $n+1$ and $n+2$.
The map is an isomorphism onto the component of $\cB^{\nu}$ except in case (1) for $g_1=g_2$ and $n=0$ and case (2),
when it is \'{e}tale of degree $2$.

For $\om_{h,S}$ a moduli stack of pointed stable curves as above, let $\pi \colon \cU_{h,S} \rightarrow \om_{h,S}$ 
denote 
the universal family, and $x_i \colon \om_{h,S} \rightarrow \cU_{h,S}$, $i \in S$, the tautological sections of $\pi$.
Define $\psi_i=x_i^*\omega_{\pi}$, the pullback of the relative dualising sheaf of $\pi$ along the section $x_i$. 
The following result is well-known, see for example \cite[Prop.~3.32]{HM98}.
\begin{lem}\label{N}
The pullback of $\cN^{\vee}$ to $\om_{g_1,S_1 \cup \{n+1\}} \times \om_{g_2,S_2 \cup \{n+2\}}$ is identified with 
$\pr_1^*\psi_{n+1} \otimes \pr_2^*\psi_{n+2}$. 
Similiarly, the pullback of $\cN^{\vee}$ to $\om_{g-1,n+2}$ is identified with $\psi_{n+1} \otimes \psi_{n+2}$. 
\end{lem}

There is an isomorphism of stacks $c : \om_{g,n+1} \rightarrow \cU_{g,n}$ 
which identifies the morphism $p_{n+1} : \om_{g,n+1} \rightarrow \omgn$
given by forgetting the last point with the projection $\pi \colon \cU_{g,n} \rightarrow \omgn$ \cite[Sec.~1--2]{Knudsen83a}.

\begin{lem}\cite[Thm.~4.1(d), p.~202]{Knudsen83b} \label{psi}
The line bundle $\psi_{n+1}$ on $\om_{g,n+1}$ is identified with the pullback of the line bundle $\omega_{\pi}(\Sigma)$
under the isomorphism $c \colon \om_{g,n+1} \rightarrow \cU_{g,n}$.
\end{lem}

\begin{cor}\label{Nbignef}
The $\bQ$-line bundle on the coarse moduli space of $\cB^{\nu}$ defined by $\cN^{\vee}$ is big and nef on each component
\end{cor}
\begin{proof}
This follows immediately from Lem.~\ref{N}, Lem.~\ref{psi}, and Thm.~\ref{bignef}. 
\end{proof}
We deduce that $H^i(\cN)=0$ for $i < \dim \cB$ by Thm~\ref{KVstack}. 
This completes the proof of Theorem~\ref{mainthm}.

\section{Proof of Theorem~\ref{coarsedefns}}

We first prove a basic result which relates the deformations of a smooth Deligne--Mumford stack and its coarse moduli space.

Let $\cX$ be a smooth proper Deligne--Mumford stack, $X$ the coarse moduli space of $\cX$, and $p : \cX \rightarrow X$ the projection.
Let $T_{\cX}$ denote the tangent sheaf of $\cX$.
Let $D \subset X$ be the union of the codimension one components of the branch locus of $p : \cX \rightarrow X$ (with its reduced structure).
Let $T_X(-\log D)$ denote the subsheaf of the tangent sheaf $T_X$ consisting of derivations which preserve the ideal sheaf of $D$.
It is the sheaf of first order infinitesimal automorphisms of the pair $(X,D)$.

\begin{lem} \label{pushfwdT}
$p_*T_{\cX}=T_X(-\log D)$
\end{lem}
\begin{proof}
The sheaves $p_*T_{\cX}$ and $T_X(-\log D)$ satisfy Serre's $S_2$ condition, and are identified over the locus where $p$ is \'{e}tale.
So it suffices to work in codimension $1$.
We reduce to the case $\cX =[\bA^1_x/\mu_e]$, where $\mu_e \ni \zeta : x \mapsto \zeta x$.
Then $X=\bA^1_x/\mu_e =\bA^1_y$, where $y=x^e$, and $D=(y=0) \subset X$.
Let $\pi : \bA^1_x \rightarrow \bA^1_x/\mu_e$ be the quotient map.
We compute 
$$p_*T_{\cX}=\left(\pi_*\cO_{\bA^1_x} \cdot \frac{\partial}{\partial x}\right)^{\mu_e}=\cO_{\bA^1_y}\cdot x\frac{\partial}{\partial x} 
=\cO_{\bA^1_y}\cdot y\frac{\partial}{\partial y} = T_X(-\log D),$$
as required.
\end{proof}

\begin{prop} \label{stack-cmod}
The first order deformations of the stack $\cX$ are identified with the first order locally trivial deformations of the pair $(X,D)$.
\end{prop}
\begin{proof}
By the Lemma, $H^1(T_{\cX})=H^1(p_*T_{\cX})=H^1(T_X(-\log D))$.
\end{proof}

We now apply this result to relate deformations of the stack $\omgn$ and its coarse moduli space $\oMgn$.

A stable $n$-pointed curve of genus $0$ has no non-trivial automorphisms.
Hence the stack $\om_{0,n}$ is equal to its coarse moduli space $\oM_{0,n}$, and $\oM_{0,n}$ is rigid by Thm.~\ref{mainthm}.
Also, recall that $\oM_{1,1}$ is isomorphic to $\bP^1$ and therefore rigid.
So, in the following, we assume that $g \neq 0$ and $(g,n) \neq (1,1)$.

Let $\sD \subset \omgn$ be the component of the boundary whose general point is a curve with two components of genus $1$ and $g-1$ 
meeting in a node, with each of the $n$ marked points on the component of genus $g-1$.
Note that each point of $\sD$ has a non-trivial automorphism given by the involution of the component of genus $1$ fixing the node.
Let $p : \omgn \rightarrow \oMgn$ be the projection, and $D \subset \oMgn$ the coarse moduli space of $\sD$.

\begin{lem}\label{aut}\cite[\S 2]{HMu82}
If $g+n \ge 4$ then the automorphism group of a general point of $\omgn$ is trivial,
and the divisor $D \subset \oMgn$ is the unique codimension $1$ component of the branch locus of $p$.
\end{lem}

Assume $g+n \ge 4$.
Let $\nu : \sD ^{\nu} \rightarrow \sD$ denote the normalisation of $\sD$, so $\sD^{\nu} = \om_{1,1} \times \om_{g-1,n+1}$.
Let $\cN_{D}$ denote the normal bundle of the map $\sD^{\nu} \rightarrow \omgn$.

\begin{lem}\label{coarseseq}
There is an exact sequence
$$0 \rightarrow T_{\oMgn}(-\log D) \rightarrow T_{\oMgn} \rightarrow p_*\nu_*\cN_{D}^{\otimes 2} \rightarrow 0.$$
\end{lem}
\begin{proof}
This is a straightforward calculation similar to \cite[Lemma, p.~52]{HMu82}.
\end{proof}

We have $H^1(T_{\oMgn}(-\log D))=H^1(T_{\omgn})=0$ by Prop.~\ref{stack-cmod} and Thm.~\ref{mainthm}.
Also $H^1(\cN_{D}^{\otimes 2})=0$ by Thm.~\ref{KVstack} because the $\bQ$-line bundle 
defined by $\cN_{D}^{\vee}$ on the coarse moduli space of $\sD^{\nu}$ is big and nef by Cor.~\ref{Nbignef}. 
So $H^1(T_{\oMgn})=0$ by Lem.~\ref{coarseseq}, that is, $\oMgn$ has no locally trivial deformations.
This concludes the proof of Thm.~\ref{coarsedefns}.

\section{Nonrigidity of moduli of surfaces} \label{example}

We exhibit a moduli space of surfaces with boundary that is not rigid.

Let $P_1,\ldots,P_4$ be $4$ points in linear general position in $\bP^2$.
Let $l_{ij}$ be the line through $P_i$ and $P_j$.
Let $l$ be a line through the point $Q=l_{12} \cap l_{34}$ such that $l$ does not pass
through $l_{13}\cap l_{24}$ or $l_{14} \cap l_{23}$ and is not equal to $l_{12}$ or $l_{34}$.
Let $S \rightarrow \bP^2$ be the blowup of the points $P_1,\ldots,P_4,Q$, and $B$ the sum of 
the strict transforms of $l$ and the $l_{ij}$ and the exceptional curves.
Then $(S,B)$ is a smooth surface with normal crossing boundary such that $K_S+B$ is very ample.
We fix an ordering $B_1,\ldots,B_{12}$ of the components of $B$.
The moduli stack $\cM$ of deformations of $(S,B)$ is isomorphic to $\bP^1\setminus\{q_1,\ldots,q_4\}$ where the $q_i$ are distinct points. Indeed, it suffices to observe that all deformations of 
$(S,B)$ are obtained by the construction above.
The moduli space $\cM$ has a modular compactification $(\ocM,\partial\ocM)$, the Koll\'ar--Shepherd-Barron--Alexeev moduli stack of stable surfaces with boundary, which is isomorphic
to $(\bP^1, \sum q_i)$.
In particular, the pair $(\ocM,\partial\ocM)$ has non-trivial deformations. 

\begin{rem}
The compact moduli space $\ocM$ is an instance of the compactifications of moduli spaces of hyperplane arrangements described in \cite{Lafforgue03} (cf. \cite{HKT06}).
\end{rem}

\section{Appendix: Kodaira vanishing for stacks}\label{app}

Let $\cX$ be a smooth proper Deligne--Mumford stack, $X$ the coarse moduli space of $\cX$, and $p : \cX \rightarrow X$ the projection.
\'{E}tale locally on $X$, $p: \cX \rightarrow X$ is of the form $p : [U/G] \rightarrow U/G$, where 
$U$ is a smooth affine variety and $G$ is a finite group acting on $U$ \cite[Lemma~2.2.3, p.~32]{AV02}.
A sheaf $\cF$ on $[U/G]$ corresponds to a $G$-equivariant sheaf $\cF_U$ on $U$, and $p_* \cF = (\pi_* \cF_U)^G$,
where $\pi : U \rightarrow U/G$ is the quotient map.

Let $\sL$ be a line bundle on $\cX$. Let $n \in \bN$ be sufficiently divisible so that for each open patch $[U/G]$ of $\cX$
as above and point $q \in U$ the stabilizer $G_q$ of $q$ acts trivially on the fibre of $\sL_U^{\otimes n}$ over $q$.
Then the pushforward $p_*(\sL^{\otimes n})$ is a line bundle on $X$.
We define $p^{\bQ}_* \sL=\frac{1}{n}p_*(\sL^{\otimes n}) \in \Pic(X) \otimes \bQ,$ and call $p^{\bQ}_*\sL$
the $\bQ$-line bundle on $X$ defined by $\sL$.

\begin{thm}\label{KVstack}
Assume that the coarse moduli space $X$ is an algebraic variety.
If the $\bQ$-line bundle $p^{\bQ}_*\sL$ on $X$ is big and nef then $H^i(\sL^{\vee})=0$ for $i < \dim \cX$.
\end{thm}

\begin{rem} 
If the coarse moduli space $X$ is smooth then Thm.~\ref{KVstack} follows from \cite[Thm.~2.1]{MO05}.
\end{rem}

Theorem~\ref{KVstack} is proved by reducing to the following generalisation of the Kodaira vanishing theorem.

\begin{thm}\cite[Thm.~2.70, p.~73]{KM98} \label{KV}
Let $X$ be a proper normal variety and $\Delta$ a $\bQ$-divisor on $X$ such that the pair $(X,\Delta)$ is Kawamata log terminal (klt).
Let $N$ be a $\bQ$-Cartier Weil divisor on $X$ such that $N \equiv M+\Delta$, where $M$ is a big and nef $\bQ$-Cartier $\bQ$-divisor.
Then $H^i(X,\cO_X(-N))=0$ for $i < \dim X$.
\end{thm}

\begin{proof}[Proof of Thm.~\ref{KVstack}]
Observe first that $X$ is a normal variety with quotient singularities.
Consider the sheaf $p_*(\sL^{\vee})$ on $X$.
If the automorphism group of a general point of $\cX$ acts nontrivially on $\sL$, then $p_*\sL^{\vee}=0$,
and so $H^i(\sL^{\vee})=H^i(p_*\sL^{\vee})=0$ for each $i$.
Suppose now that the automorphism group of a general point acts trivially on $\sL$.
Then $p_* \sL^{\vee}$ is a rank 1 reflexive sheaf on $X$. Write $p_*\sL^{\vee}=\cO_X(-N)$, where $N$ is a Weil divisor on $X$.
Let $n \in \bN$ be sufficiently divisible so that $p_*^{\bQ}(\sL)=\frac{1}{n}p_*(\sL^{\otimes n})$ as above.
Let $M$ be a $\bQ$-divisor corresponding to the $\bQ$-line bundle $p_*^{\bQ}\sL$.
There is a natural map $(p_* \sL^{\vee})^{\otimes n} \rightarrow p_*(\sL^{\vee \otimes n})$, i.e., a map
$\cO_X(-nN) \rightarrow \cO_X(-nM)$, which is an isomorphism over the locus where $p$ is \'{e}tale.
So $N \equiv M + \Delta$, where $\Delta$ is an effective $\bQ$-divisor supported on the branch locus of $p$.
Let $D_1,\ldots,D_r$ be the codimension $1$ components of the branch locus.
Let $e_i$ be the ramification index at $D_i$, and $a_i$ the age of the line bundle $\sL^{\vee}$ along $D_i$.
That is, after removing the automorphism group of a general point of $\cX$, a transverse slice of $\cX$ at a general point of $D_i$
is of the form $[\bA^1_x / \mu_{e_i}]$, where $\mu_{e_i} \ni \zeta : x \mapsto \zeta \cdot x$,
and $\mu_{e_i}$ acts on the fibre of $\sL^{\vee}$ by the character $\zeta \mapsto \zeta^{-a_i}$, where $0 \le a_i \le e_i-1$.
We compute that $\Delta = \sum \frac{a_i}{e_i}D_i$.

We claim that $(X,\Delta)$ is klt.
Let $\Delta'=\sum \frac{e_i-1}{e_i}D_i$, then $K_{\cX}=p^*(K_X+\Delta')$, and $\cX$ is smooth, so $(X,\Delta')$ is klt by 
\cite[Prop.~5.20(4), p.~160]{KM98}. Now $\Delta \le \Delta'$ and $X$ is $\bQ$-factorial, so $(X,\Delta)$ is also klt.
We deduce that $H^i(\sL^{\vee})=H^i(p_*\sL^{\vee})=H^i(\cO_X(-N))=0$ for $i< \dim \cX$ by Thm.~\ref{KV}.
\end{proof}

\medskip
\noindent
Paul Hacking,  Department of Mathematics, University of Washington, Box 354350, Seattle, WA~98195; 
\texttt{hacking@math.washington.edu} \\

\end{document}